\documentclass[10pt,a4paper,oneside,headinclude,footinclude,BCOR5mm]{article}

\usepackage{fancyhdr,amssymb,amsbsy,amsthm}

\title{\normalfont{The general definition of K\"ahler angle}}

\author{Yongpin Zhu\textsuperscript{1}}

\date{}

\begin{document}


\renewcommand{\sectionmark}[1]{\markright{#1}} 

\pagestyle{fancy} 


\maketitle 






\section*{Abstract}

In this paper, we give a general definition of K\"ahler angle. There are many results about K\"ahler angle one can try to generalize to the general case. 



\let\thefootnote\relax\footnotetext{\textsuperscript{1} \textit{School of Mathematical Sciences, University of Science and Technology of China, Hefei 230026, P.R.of China}}


\section{Introduction}

 It is well-konwn that the algebraic structure of a vector space can be generalized to the corresponding objects on the manifold by a pointwise way. The geometric version of the algebraic property of Wirtinger's inequality was given, such as, in \cite{L}. In this paper, we give the general definition of the K\"ahler angle by the geometric version of Wirtinger's inequality.

Suppose that $(N,J,\langle,\rangle)$ is a Hermite manifold, where $\langle,\rangle$ is the Riemannnian metric on $N$, $J$ is a complex structure compatible with $\langle,\rangle$, i.e.$$\langle JU,JV\rangle=\langle U,V\rangle.$$ The K\"ahler form $\omega$ on $N$ is defined by $$\omega(U,V)=\langle JU,V\rangle.$$



\section{The general definition of K\"ahler angle}

For the convenience of discussion, we first state a algebraic version of Wirtinger's inequality which takes from \cite{H}, the proof is almost the same as the geometric version.

\newtheorem{theorem}{Theorem}[section] 

\begin{theorem}
Let $(V,J,\langle,\rangle)$ be a euclidian vector space endowed with a compatible almost complex structure $J$, and the associated fundamental form $\omega$. Let $W\subset V$ be an oriented subspace of dimension $2m$. The induced scalar product on $W$ together with the chosen orientation define a natural volume form $vol_{W}\in \bigwedge^{2m}W^{*}$. Then  $$\omega^ {m}|_{W}\le{m!Vol_{W}},$$
and that equality holds if and only if $W\subset V$ is a complex subspace, i.e.$J(W)=W$, and the orientation is the one induced by the almost complex structure.
\end{theorem}

The geometric version of Wirtinger's inequality was given by a pointwise way. 

\begin{theorem}[Wirtinger inequation]
  Let $N$ be a $n$-dimensional Hermite manifold, let $M\subset{N}$ is any $2m$-dimensional orientable real submanifold. At any point $p\in{M}$, let $dV_{p}$ denote the volume form of the induced metric on $M$. Then the restriction of the $m$-power $\omega^{m}=\omega\wedge\cdots\wedge\omega$ of the K\"ahler form of $N$ to $T_{p}M$ satisfies 
  $$\frac{\omega^{m}}{m!}\le{dV_{p}},$$
  and equality holds if and only if $T_{p}M$ is a complex subspace of $T_{p}N$, with the canonical orientation.
\end{theorem}

\begin{proof}
If $m=1$, then for any unit vectors $v,w\in T_{p}N$,
$$\omega(w,v)^{2}=\langle Jv,w\rangle^{2}\le |Jv|^{2}|w|^{2}=|v|^{2}|w|^{2}=1.$$
and that equality holds $\iff$ $w=Jv$, i.e. $span\{w,v\}$ is the one-dimensional complex subspace of $T_{p}N$.
For the general case, we now consider the form $\omega$ restircted to $T_{p}M$, if we choose the oriented orthonormal basis $\{e_{1},\cdots,e_{2m}\}$ of $T_{p}M$, then $\omega(e_ {i},e_{j})=\langle Je_{i},e_{j}\rangle=J^{j}_{i}$. By $J^{2}=-id$ and the elementary linear algebra, we know that there exists the oriented orthonormal basis, which still denote by $\{e_{1},\cdots,e_{2m}\}$, such that $J$ can be represented by 
$$\left(
\begin{array}{ccccc}
   0 & \lambda_{1} &  &  &  \\
    -\lambda_{1} & 0 &  &  &  \\
    &  & \ddots&  &   \\
    &  &  & 0 & \lambda_{m} \\
    &  &  & -\lambda_{m} & 0
\end{array}\right),
$$
Then we have $\omega(e_{2k-1},e_{2k})=\lambda_{k}, k=1,2,\cdots,m$. Let $\{\theta^{i}\}$ be the dual basis of $\{e_{i}\}$, then
$$\omega=\sum\limits_ {k=1}^{m}\lambda_{k}\theta^{2k-1}\wedge\theta^{2k},$$ 
and $|\lambda_{k}|\le 1$. Thus, we have
$$\omega^{m}=m!\lambda_{1}\cdots\lambda_{m}\theta^{1}\wedge\cdots\wedge\theta^{2m}=m!\lambda_{1}\cdots\lambda_{m}dV_{p},$$ 
so $|\omega^{m}|\le m!dV_{p}$, and equality holds $\iff$ $|\lambda_{1}\cdots\lambda_{m}|=1\iff \omega(e_{2k-1},e_{2k})=\pm{1},\forall{k=1,\cdots,m}\iff e_{2k-1}=\pm Je_{2k}, 1\le k\le m \iff T_{p}M$ is a complex subspace of $T_{p}N$. Equality holds without the absolute value sign exactly when the orientation agree with the one induced by $J$.  
\end{proof}

Note that the statement of the theorem is pointwise, so it still holds true for the almost-Hermite manifolds, in fact, it is a geometric version of the algebraic version of Wirtinger's inequality. As the top form over $M$, note that the dimension of the top form space is $1$, we can introduce the function 
 $$\frac{\omega^ {m}|_{T_{p}M}}{m!dV_{p}}\le{1}$$ 
on $M$ to characterize the difference between the two top form. When $n=2, m=1$, the function is $\cos{\alpha}$ in \cite{CW}, which was first introduced in \cite{CW}. And then $\alpha$ was called the K\" ahler angle because of its geometric meaning in that case: $\cos{\alpha}=\omega(e_{1},e_{2})=\langle Je_{1},e_{2}\rangle$, i.e. the angle between the unit vector $Je_{1}$ and $e_{2}$, where $e_{1}$, $e_{2}$ are local orthonormal frame fields on $M$. Therefore, we still denote the function by $\cos{\alpha}$ and call $\alpha$ the K\"ahler angle.

\newtheorem{definition}[theorem]{Definition}

\begin{definition}[K\"ahler angle]
  Let $N$ be a $n$-dimensional almost-Hermite manifold, let $M\subset{N}$ is any $2m$-dimensional orientable real submanifold. Then we define K\"ahler angle $\alpha$ by 
  $$\cos{\alpha} = \frac{\omega^ {m}|_{TM}}{m!dV}, \eqno(*)$$
  where $\omega$ is the K\"ahler form of $N$, $dV$ denote the volume form of the induced metric on $M$.
\end{definition}

Next, we consider the case that $\cos{\alpha}\equiv{1}$. First, we recall the well-known $Newlander-Nirenberg$ theorem.

\begin{theorem}[Newlander-Nirenberg]
  Let $(M,J)$ be an almost-Hermite manifold, if $J$ is integrabel, i.e.the torsion tensor of $J$ is $0$, then $J$ is a complex structure induced by a complex manifold structure of $M$. 
\end{theorem}

When $\cos{\alpha}\equiv{1}$, by the Wirtinger inequality, for any $p\in{M}$, $T_{p}M$ is a complex subspace of $T_{p}N$, i.e. $J_{p}(T_{p}M)\subset{T_{p}M}$, thus $(M, J|_{M}, g)$ is a almost-Hermite submanifold, where $g$ is an induced metric of $M$. If $N$ is a Hermite manifold, then $(M, J|_{M}, g)$ is a Hermite submanifold. In fact, by the definition of the torsion tensor of the almost complex structure, the torsion tensor of $J|_{M}$ is the limit of the torsion tensor of the complex structure $J$ on $M$, so it is $0$. Then it is known from the $Newlander-Nirenberg$ theorem that it is a complex submanifold. We call the point of $\alpha=0$ (i.e. $\cos{\alpha}=1$) a complex point, and the point of $\alpha=\pi$ (i.e. $\cos{\alpha}=-1$) a anti-complex point.

By the definition, K\"ahler angles $\alpha\in[0,\pi]$ are determined by $(*)$ and the function $\arccos$, thus the function $\alpha$ is continuous on $M$ and smooth except for the complex points and the anti-complex points. When $m=1$, and we consider it locally, we can always allow the angle of the two vector fields to take a value outside $[0,\pi]$, and when we calculate the derivative, the two angle differ by at most one negative sign, so that $|\nabla\alpha|$ is consistent as a function on $M$. And when we consider the subsequent question, we always get the derivative of angle $\alpha$ by the derivative of the function $\cos{\alpha}$, so we can also directly consider the function $(*)$. But at least, when $m=1$, expressing the function $(*)$ as $\cos{\alpha}$ can effectively simplify the expression by using the properties of $\cos{\alpha}$ (e.g., $|\cos(\alpha)|\le{1}$, $\cos{\alpha}^{2}+\sin{\alpha}^{2}=1$, etc.).
 
In the general dimensionality, after a choice of the appropriate local unitary frames, the function$(*)=\cos{\alpha_{1}}\cdots\cos{\alpha_{m}}$, where $\cos{\alpha_{k}}=\omega(e_{2k-1},e_{2k}), k=1,\cdots,m$. Thus the $\alpha_{k}$ can be treated similarly, note that $\alpha_{k}$ depend on the choice of the local unitary frames, but not for K\"ahler angle $\alpha$, so we study the function $(*)$ directly. For convenience, we call the function $(*)$ is the K\"ahler function.

In this paper, we only focus on the definition of K\"ahler angle. But after giving the general definition of K\"ahler angle, the original results about K\"ahler angle can be generalized from $N$ is a K\"ahler manifold , and $m=1, n=2$, to the general case: $N$ is a almost-Hermite manifold, and any dimension $2m$ and any codimension $2(n-m)$.



\newpage

\begin{thebibliography}{99}



\bibitem{CW} S. S. Chern and J. Wolfson, Minimal surfaces by moving frames, Amer. J. Math. 105(1983), 59-83.





\bibitem{L} Lawson, H.B. Lectures on Minimal Submanifolds. 1. University of California, 1980.

\bibitem{H} D.Huybrechts, Complex geometry. An introduction [M]. J.H.w. Dietz nachf. g.m.b.h, 2005.


\end{thebibliography}
\end{document}